\documentclass{amsart}

\usepackage{amsmath}
\usepackage{amsfonts}
\usepackage{amssymb}
\usepackage{amsthm}
\usepackage{comment}
\usepackage{epsfig}
\usepackage{psfrag}
\usepackage{mathrsfs}
\usepackage{amscd}
\usepackage[all]{xy}
\usepackage{rotating}
\usepackage{lscape}
\usepackage{amsbsy}
\usepackage{verbatim}
\usepackage{moreverb}
\usepackage{xspace}

\newtheorem{theorem}{Theorem}[section]
\newtheorem{prop}[theorem]{Proposition}
\newtheorem{lemma}[theorem]{Lemma}

\theoremstyle{definition}

\theoremstyle{remark}

\newcommand{\nn}{\nonumber}
\newcommand{\nid}{\noindent}
\newcommand{\ra}{\rightarrow}

\newcommand{\xra}{\xrightarrow}
\newcommand{\xla}{\xleftarrow}

\DeclareMathOperator*{\colim}{colim}

\DeclareMathOperator{\Map}{Map}
\DeclareMathOperator{\coker}{coker}

\newcommand{\hra}{\hookrightarrow}
\newcommand{\sm}{\wedge}

\newcommand{\tw}{\mathrm{tw}}
\newcommand{\Spec}{\mathrm{Spec}}
\newcommand{\fU}{\mathfrak{U}}
\newcommand{\fg}{\mathfrak{g}}
\newcommand{\fh}{\mathfrak{h}}
\newcommand{\ft}{\mathfrak{t}}
\newcommand{\aff}{\mathrm{aff}}
\newcommand{\Spinc}{{Spin^c}}

\def\llarrow{   \hspace{.05cm}\mbox{\,\put(0,-2){$\leftarrow$}\put(0,2){$\leftarrow$}\hspace{.45cm}}}

\def\lllarrow{  \hspace{.05cm}\mbox{\,\put(0,-3){$\leftarrow$}\put(0,1){$\leftarrow$}\put(0,5){$\leftarrow$}\hspace{.45cm}}}

\def\llllarrow{  \hspace{.05cm}\mbox{\,\put(0,-6){$\leftarrow$}\put(0,-2){$\leftarrow$}\put(0,2){$\leftarrow$}\put(0,6){$\leftarrow$}\hspace{.45cm}}}

\newcommand{\lla}{\llarrow}
\newcommand{\llla}{\lllarrow}
\newcommand{\lllla}{\llllarrow}
\newcommand{\dra}{\Rightarrow}

\def\cA{\mathcal A}
\def\cH{\mathcal H}
\def\cK{\mathcal K}

\def\CC{\mathbb C}

\def\NN{\mathbb N}

\def\ZZ{\mathbb Z}

\begin{document}

\author{Christopher L. Douglas} 
\thanks{The author was supported in part by an NSF Graduate Research Fellowship and in part by a Miller Research Fellowship.}
\address{Department of Mathematics, University of California, Berkeley, CA 94720, USA}
\email{cdouglas@math.berkeley.edu}

\title{On the Structure of the Fusion Ideal}

\begin{abstract}

We prove that there is a finite level-independent bound on the number of relations defining the fusion ring of positive energy representations of the loop group of a simple, simply connected Lie group.  As an illustration, we compute the fusion ring of $G_2$ at all levels.

\end{abstract}

\maketitle

%
%
%
%
%

\vspace*{-15pt}
\section{Introduction} \label{sec-introduction}

\subsubsection*{Background and overview}

Two-dimensional rational conformal field theory has applications to the study of critical phenomena in statistical mechanics and condensed matter physics, and is intimately related to the theory of operator algebras and to three-dimensional topological field theory~\cite{bpz, henkel, tsvelik, wassermann, witten}.  The basic invariant of a rational conformal field theory is the associated fusion ring, whose generators are the primary fields of the theory and whose multiplication is determined by the dimensions of the spaces of conformal blocks~\cite{verlinde}.  For a class of conformal field theories associated to loop groups, the fusion ring can be described as a Grothendieck group of positive energy representations, with product structure determined by Connes fusion.  Our purpose in this note is to investigate algebraic properties of these fusion rings, and particularly to describe a certain bound on the complexity of the fusion ring of a loop group.

The fusion ring $F_k[G]$ of positive energy representations of a level $k$ central extension of the loop group $LG$ of the compact, simple, simply connected Lie group $G$ can be presented as a quotient $R[G]/I_k$ of the representation ring of $G$ by the ``fusion ideal" $I_k$.  This quotient is freely generated as an abelian group by the irreducible representations of $G$ whose highest weights have level less than or equal to $k$.  The fusion ideal contains and is often generated by the level $k+1$ irreducible representations of $G$, suggesting that the number of generators of the ideal goes to infinity with the level.  However, Gepner~\cite{gepner} conjectured that these fusion rings are all global complete intersections, that is the fusion ideal is generated by exactly $n$ elements, for a group $G$ of rank $n$.  We describe a situation intermediate between these two extremes: there is a finite level-independent upper bound on the number of representations needed to generate the fusion ideal.  Moreover, there is a uniform presentation of the fusion ring of $G$ in terms of bases for the modules of representations of central extensions of the centralizer subgroups of $G$.  As an example, we work out the necessary bases for the subgroups of $G_2$ and thereby give an explicit computation of the $G_2$ fusion ring.

Our primary technique is to analyze a twisted Mayer-Vietoris spectral sequence converging to the twisted equivariant $K$-homology of the group $G$.  By results of Freed, Hopkins, and Teleman~\cite{fht1, fht2, fht3} the twisted equivariant $K$-homology $K^G_{k+h^{\scriptscriptstyle\vee}}(G)$ for twisting $k+h^{\scriptscriptstyle\vee}$ is isomorphic to the fusion ring $F_k[G]$ of positive energy representations of the loop group $LG$ at level $k$.  In fact, the twisted equivariant $K$-homology of $G$ is a Frobenius algebra encoding the full two-dimensional topological field theory associated to the underlying conformal field theory.  Moreover this twisted equivariant $K$-homology admits operations by $K$-theory classes of moduli spaces of surfaces, altogether forming the state space of a topological conformal field theory~\cite{fht-orient}.  Though here we limit our attention to the multiplicative behavior of the fusion ring, computational tools from twisted homotopy theory could also be used to analyze the further topological conformal structure of these field theories.

\subsubsection*{Results and organization}

Parametrized stable homotopy theory is a convenient framework for studying twisted generalized cohomology theories, and in section~\ref{sec-teht} we recall the definitions of parametrized spectra and their associated twisted homology and cohomology theories.  We then introduce the twisted Mayer-Vietoris spectral sequence in section~\ref{sec-twmvss} and discuss a particular such spectral sequence  for the $K$-homology and equivariant $K$-homology of $G$.  In section~\ref{sec-resolution} we reformulate the equivariant spectral sequence in terms of ``twisted representation modules", that is in terms of modules of representations of central extensions of subgroups of $G$, and we express the $d^1$ differentials of the spectral sequence as twisted holomorphic induction maps between these modules.  The resulting resolution of the fusion ring is implicit in the work of Freed, Hopkins, and Teleman, is described in papers of Kitchloo and Morava~\cite{kitchloomorava, kitchloo}, and is presented in detail from the point of view of the twisted $K$-theory of $C^*$-algebras by Meinrenken~\cite{meinrenken}.

In section~\ref{sec-g2}, we compute bases for the representation modules of the centralizer subgroups of $G_2$, with particular attention to the twisted module of representations of the nontrivial central extension of the $SO(4)$ subgroup.  We then exploit singular planes for the twisted holomorphic induction maps to simplify the analysis of the differential in the spectral sequence, thereby facilitating an explicit description of the fusion ring.
\begin{theorem}
For $k>0$ the level $k$ fusion ring of $G_2$ is presented as follows:
\begin{equation} \nn
F_k[G_2] =
\begin{cases}
R[G_2]/(\rho_{(1, \frac{k}{2})}; \rho_{(0, \frac{k}{2})}+\rho_{(0,\frac{k}{2} + 1)}; \rho_{(1,\frac{k}{2} - 1)}+\rho_{(1,\frac{k}{2} +1)}; \rho_{(k+2,0)})
& \text{$k$ even} \\
R[G_2]/(\rho_{(0, \frac{k+1}{2})};\rho_{(1, \frac{k-1}{2})}+\rho_{(1,\frac{k+1}{2})}; \rho_{(0, \frac{k-1}{2})}+\rho_{(0,\frac{k+3}{2})}; \rho_{(k+2,0)})
& \text{$k$ odd}
\end{cases}
\end{equation}
Here $\rho_{(a,b)}$ is the irreducible representation with highest weight $\lambda_s^a \lambda_l^b$, for $\lambda_s$ and $\lambda_l$ the short and long fundamental weights of $G_2$.
\end{theorem}

We abstract this $G_2$ computation to other simply-connected groups in section~\ref{sec-finiteness}.  In particular, we establish that the $E^1$ term of the twisted Mayer-Vietoris spectral sequence for the twisted equivariant $K$-homology of $G$ consists entirely of free $R[G]$-modules.  This entails a generalization of a result of Pittie~\cite{pittie}, which might be of independent interest:
\begin{prop}
Let $G$ be simple and simply connected, and let $H$ be the centralizer of an element of $G$.  The restriction of the generator of $H^3_G(G;\ZZ)$ to the orbit $G/H$ classifies an $S^1$-central extension of $H$, denoted $\widetilde{H}$.  The submodule $R_k[H] \subset R[\widetilde{H}]$, consisting of representations of $\widetilde{H}$ for which the central circle acts by the $k$-th power of scalar multiplication, is free as an $R[G]$-module.
\end{prop}
For convenient reference, we list all the centralizer subgroups $H$ of all simple, simply connected groups $G$ for which the twisted representation module $R_k[H]$ is in fact twisted, that is not isomorphic to the representation ring $R[H]$.  The fusion ring admits a presentation in terms of induction maps on bases for these twisted modules, and the level-independent finiteness of the fusion ideal follows.
\begin{theorem}
For $G$ a compact, simple, and simply connected Lie group, there exists an integer $n_G$ such that for all levels $k$, the fusion ring of $G$ at level $k$ has a presentation of the form $F_k[G] = R[G]/(\rho_1(k), \rho_2(k), \ldots, \rho_{n_G}(k))$, for representations $\rho_i(k) \in R[G]$ depending on the level.  An upper bound on the number of generators $n_G$ of the fusion ideal is $\sum_{i \in \bar{D}} |W_G|/|W_{Z(\Lambda_i)}|$; here we denote the Dynkin diagram of $G$ by $\bar{D}$, the centralizer of the edge of the Weyl chamber containing the fundamental weight $\lambda_i$ by $Z(\Lambda_i)$, and the Weyl group of $H \subset G$ by $W_H$.
\end{theorem}

We have recently discovered two papers that also study the algebraic structure of the fusion ring: the first~\cite{bouwridout}, by Bouwknegt and Ridout, shows that in type $A$ and type $C$ the fusion ring is always a global complete intersection over the integers---in fact the authors construct an integral fusion potential in these cases; the second paper~\cite{boysalkumar}, by Boysal and Kumar, describes a series of intriguing conjectures for finite but not complete intersection presentations of the fusion rings for the classical groups and for $G_2$.

\subsubsection*{Acknowledgments}

We would like to thank Mike Hopkins for introducing us to twisted $K$-theory, and for suggesting the twisted equivariant Mayer-Vietoris spectral sequence as an approach, initially, to computations in non-equivariant twisted $K$-theory.  We became distracted by other techniques for those computations, and we apologize for the resulting extreme delay in the publication of these results.  We thank Eckhard Meinrenken for various helpful conversations and particularly for sharing his $C^*$-algebra perspective on the twisted Mayer-Vietoris spectral sequence.  We are also indebted to Mark Haiman, Max Lieblich, and Andre Henriques for informative and illuminating discussions.

\section{Computing twisted equivariant $K$-theory} \label{sec-computingtekt}

We begin this section by reviewing a general framework for twisted equivariant homology theory, namely parametrized equivariant spectra.  We then introduce our primary computational tool, the twisted Mayer-Vietoris spectral sequence, and describe it in detail for the elementary twists of equivariant $K$-theory over a simple simply connected compact Lie group.  Finally, we reformulate this spectral sequence in representation-theoretic terms, obtaining a resolution of the fusion ring by explicit twisted representation modules.

\subsection{Recollections on twisted equivariant homology theories} \label{sec-teht}

A function on a space $X$ is a map $X \ra \CC$, and a twisted function is a section $\xymatrix@1{\CC \ar[r] & L \ar@<-.5ex>[r] & X \ar@<-.5ex>@{-->}[l]}$ of a line bundle over $X$.  Similarly, for a spectrum $F$, an $F$-cohomology class is a homotopy class of maps $X \ra F$, and a twisted $F$-cohomology class is a homotopy class of sections $\xymatrix@1{F \ar[r] & E \ar@<-.5ex>[r] & X \ar@<-.5ex>@{-->}[l]}$ of a ``bundle of spectra" over $X$ with fiber $F$.  If $F$ is instead a $G$-equivariant spectrum, for some compact Lie group $G$, the same heuristic suggests a notion of twisted $G$-equivariant $F$-cohomology, as equivariant sections of a bundle of equivariant spectra.

These ideas can be made precise using parametrized equivariant spectra, a theory of which is developed in extensive detail by May and Sigurdsson~\cite{maysigurdsson}.  Recall that a prespectrum is a collection of based spaces $\{F_i\}$ together with based maps $\Sigma F_i \ra F_{i+1}$.  A prespectrum parametrized over $X$ is a collection of spaces $\{E_i\}$ equipped with maps $E_i \ra X$ and sections $X \ra E_i$, together with structure maps $\Sigma_X E_i \ra E_{i+1}$ commuting with the maps to and from $X$.  The cohomology and homology of $X$ with coefficients in $F$ and $E$ are defined as follows:
\begin{align}
F^n(X) &:= \pi_{-n} \Map(X,F) = \colim \pi_{i-n} \Map(X,F_i)    \nn\\ 
E^n(X) &:= \pi_{-n} \, \Gamma(X, E) = \colim \pi_{i-n} \, \Gamma(X,E_i)    \nn\\
F_n(X) &:= \pi_{n} \, X_+ \sm F = \colim \pi_{i+n} \, (X \times F_i)/X    \nn\\
E_n(X) &:= \pi_{n} \, E/X = \colim \pi_{i+n} \, E_i/X    \nn
\end{align} 
\nid Notice that $E$ is not itself a spectrum, so the notation cannot cause confusion: $E^n(X)$ and $E_n(X)$ are twisted forms of $F$-cohomology and $F$-homology of $X$.  Here in the first instance the expressions $\Map$, $\Gamma$, $\sm$, and $/X$ refer respectively to the derived mapping space, sections, smash product, and total prespectrum.  The nonderived prespectrum of sections is defined levelwise by $\Gamma(X,E)_i = \Gamma(X,E_i)$; the nonderived total prespectrum is defined similarly by $(E/X)_i = E_i/X$.  The section and total functors are respectively the right adjoint $p_*$ and the left adjoint $p_!$ to the pullback functor $p^*: \Spec \ra \Spec/X$ for the projection $p:X \ra *$; that is $\Gamma(X,E) = p_*(E)$ and $E/X = p_!(E)$.  Under appropriate assumptions on $X$, $F$, and $E$, the above four colimit expressions, interpreted without derivation, give a more concrete description of the corresponding untwisted and twisted cohomology and homology groups.

The twisted equivariant story is analogous.  For a compact Lie group $G$, a $G$-equivariant prespectrum is a collection of based $G$-spaces $\{F_V\}$ for $G$-representations $V$, together with compatible equivariant structure maps $\Sigma^W F_V \ra F_{V \oplus W}$.  Similarly, a $G$-equivariant prespectrum parametrized over the $G$-space $X$ is a collection of $G$-spaces $\{E_V\}$ over $X$, with sections from $X$, together with compatible equivariant structure maps $\Sigma^W_X E_V \ra E_{V \oplus W}$.  The equivariant cohomology and homology of $X$ with coefficients in $F$ and $E$ are defined as follows:
\begin{align}
F^n_G(X) &:= \pi_{-n} \Map(X,F)^G  \nn\\
E^n_G(X) &:= \pi_{-n} \, \Gamma(X,E)^G  \nn\\
F_n^G(X) &:= \pi_{n} (X_+ \sm F)^G  \nn\\
E_n^G(X) &:= \pi_{n} (E/X)^G  \nn
\end{align}
These compact expressions are a convenient mnemonic, though they suppress various details.  As before, the expressions $\Map$, $\Gamma$, $\sm$, and $/X$ should all be understood as derived.  Similarly $(-)^G$ is a derived functor, specifically the derivation of the composite $G\textrm{-Spec}/X \xra{\text{forget}} \textrm{naive-}G\textrm{-Spec}/X \xra{\text{$G$-fixed}} \textrm{Spec}/X$---the first functor forgets to naive $G$-prespectra, that is $G$-prespectra indexed only on trivial $G$-representations, and the second functor takes levelwise $G$-fixed points.  Note that even if $E$ is a fibrant $G$-prespectrum over $X$, the total spectrum $E/X$ is unlikely to be fibrant, and so must be replaced by a fibrant $G$-prespectrum before taking levelwise fixed points---this complexity appears already in nonparametrized equivariant homology theory.

\subsection{The twisted Mayer-Vietoris spectral sequence} \label{sec-twmvss}

Twisted homology theories enjoy all the good properties of ordinary homology theories, including the existence of a Mayer-Vietoris spectral sequence.  Indeed the Mayer-Vietoris spectral sequence, which captures the idea that we can recover the value of a theory globally from the values locally, encodes the essence of what it means to be a homology theory.

Suppose we have a $G$-prespectrum $E$ parametrized over a $G$-space $X$.  Given an open equivariant covering $\fU=\{U_i \hra X\}$ of $X$, we can form the associated simplicial cover 
$$s\fU = \coprod_i U_i \lla \coprod_{i,j} U_{ij} \llla \coprod_{i,j,k} U_{ijk} \lllla \cdots$$ 
Here $U_I$ refers to the intersection of the $U_i$ for $i \in I$.  Note that the indexing set $I$ may contain repeated indices, and the simplicial object $s\fU$ does have degeneracy maps.  Provided the cover is numerable, the realization $|s\fU|$ is homotopy equivalent to $X$, and we may hope to obtain information about the twisted equivariant homology $E^G_*(X)$ from the homology $E^G_*(U_I)$ of the pieces of the simplicial cover.  Indeed, the simplicial filtration of $|s\fU|$ leads, a la Segal~\cite{segal-csss}, to our desired spectral sequence:

\begin{prop}
For $E$ a $G$-prespectrum parametrized over a $G$-space $X$, and $\fU$ a numerable open covering of $X$ with associated simplicial cover $s\fU$, there are spectral sequences
\begin{align}
E^2_{pq} &= H_p(E^G_q(s\fU)) \dra E^G_{p+q}(X)  \nn\\
E_2^{pq} &= H^p(E_G^q(s\fU)) \dra E_G^{p+q}(X)  \nn
\end{align}
\end{prop}

We think of the first spectral sequence as a twisted Mayer-Vietoris spectral sequence, and of the second as a twisted Bousfield-Kan spectral sequence.  May and Sigurdsson~\cite{maysigurdsson} refer to both as ``Mayer-Vietoris", but Freed, Hopkins, and Teleman~\cite{fht1} consider the second of ``Atiyah-Hirzebruch" type.  Since both spectral sequences are constructed by a method of Segal, we see that, in any case, it's a party.  

\subsubsection{TMVSS in $K$-homology} \label{sec-inkhom}

Here we describe the twisted Mayer-Vietoris spectral sequence for the twisted $K$-homology of a simple, simply connected Lie group $G$, and in the next section consider the corresponding equivariant spectral sequence.

Certain twists of $K$-theory can be constructed as follows.  Fix a Hilbert space $\cH$ and an action of the projective unitary group $PU(\cH)$ on a fixed model $K$ for the $K$-theory prespectrum; such a model can be built from spaces of Fredholm operators on $\cH$, as in Atiyah-Singer~\cite{atiyahsinger-skew} or more recently Atiyah-Segal~\cite{atiyahsegal-tkt}.  Given a principal $PU(\cH)$ bundle $P$ on a space $X$, form the parametrized prespectrum $E_P := P \times_{PU(\cH)} K$ over $X$---this parametrized prespectrum has spaces the levelwise associated bundles $P \times_{PU(\cH)} K_i$.  Such principal $PU(\cH)$ bundles over $X$ are classified by maps $X \ra BPU(\cH) \simeq K(\ZZ,3)$ and therefore up to isomorphism by classes in $H^3(X;\ZZ)$.  Given a twisting class $\tau \in H^3(X;\ZZ)$, define the $\tau$-twisted $K$-homology of $X$ as $K_{\tau,i}(X) := (E_P)_i(X)$, for a principal bundle $P$ classified by $\tau$.

Fix a simple simply connected Lie group $G$ of rank $n$.  In order to describe the twisted Mayer-Vietoris spectral sequence for $K_{\tau,n}(G)$ for a twisting $\tau \in H^3(G;\ZZ) \cong \ZZ$, we construct a particular cover $\fU = \{U_i\}_{i=0}^n$ of $G$ as follows.  We use the following notation:
\begin{align*}
T &= \text{maximal torus of } G \\
\fg &= \text{Lie algebra of } G \\
\ft &= \text{Lie algebra of } T \\
\Lambda_R &= \text{root lattice of } \fg \\
\Lambda_W &= \text{weight lattice of } \fg \\
\Gamma_R &= \text{coroot lattice of } \fg \\
\Gamma_W &= \text{coweight lattice of } \fg \\
W &= \text{Weyl group of } \fg \\
W^{\aff} &= \text{affine Weyl group of } \fg \\
A &= \text{Weyl alcove of } \fg \\
\{\alpha_i\}_{i=1}^n &= \text{simple roots of } \fg \\
\{\lambda_i\}_{i=1}^n &= \text{fundamental weights of } \fg \\
D &= \text{affine Dynkin diagram of } \fg \\
h^{\scriptscriptstyle\vee} &= \text{dual Coxeter number of } \fg  
\end{align*}
Recall that the root lattice $\Lambda_R \subset \ft^*$ is generated by the roots $\{\alpha\}$, the eigenvalues for the eigenspaces $\fg_{\alpha}$ of the restriction of the adjoint representation of $\fg$ to the maximal torus $\ft$.  The coroot lattice $\Gamma_R \subset \ft$ is generated by the coroots $\{h_{\alpha}\}$, which by definition are the unique elements $h_{\alpha} \in [\fg_{\alpha},\fg_{-\alpha}]$ such that $\alpha(h_{\alpha}) = 2$.  The weight lattice $\Lambda_W \subset \ft^*$ is by definition evaluation dual to the coroot lattice, and the coweight lattice $\Gamma_W \subset \ft^*$ is by definition evaluation dual to the root lattice.  As $G$ is simply connected, the coroot lattice can also be characterized as $\Gamma_R = \ker(\exp: \ft \ra G)$.  The roots $\alpha \in \ft^*$ of $\fg$ act on $\ft$ by reflections in the kernels $\ker \alpha$, and these reflections generate the Weyl group $W$.  The affine Weyl group is the semi-direct product $\Gamma_R \rtimes W$, and the Weyl alcove is the quotient $A = \ft/W^{\aff}$.  

The group $G$ acts on itself by conjugation, with quotient $G/G$ isomorphic to the Weyl alcove $A$---the isomorphism is the exponential map $A \xra{\exp} G/G$.  The cover of $G$ we are interested in is pulled back along $\pi: G \ra G/G$ from a cover of the alcove.  The Weyl alcove $A$ is a simplex embedded in $\ft$.  The $n+1$ vertices $\{v_{\hat{\imath}}\}_{i=0,\ldots,n}$ of $A$ are determined as follows: the vertex $v_{\hat{\imath}}$ is the fixed point of the subgroup of the affine Weyl group generated by the reflections corresponding to the roots $\alpha_0, \alpha_1, \ldots, \widehat{\alpha_i}, \ldots, \alpha_n$.  Here $\alpha_0$ is the affine (lowest) root, which acts on $\ft$ by reflection in the plane $P$ perpendicular to the coroot of the highest root, with $P$ half way between the origin and that coroot.  We have used the shorthand $\hat{\imath}$ to denote the complement in $\{0, 1, \ldots, n\}$ of the element $i$; indeed we will abbreviate the complement of any collection $S \subset \{0, 1, \ldots, n\}$ by $\widehat{S}$.  Generally, for a subset $S$ of the simple affine roots $\{\alpha_i\}_{i=0}^n$ (that is a subset of the affine Dynkin diagram $D$), there is a corresponding face $F_S$ of the alcove, defined as the fixed set of the subgroup of $W^{\aff}$ generated by the reflections corresponding to $\{\alpha_i\}_{i \in S}$.  

Cover the alcove $A$ by open sets $\{\widetilde{U_{\hat{\imath}}}\}$, with $\widetilde{U_{\hat{\imath}}} = A \backslash F_{i}$---these open sets are the complements of the codimension-one faces---and define the cover $\fU=\{U_{\hat{\imath}}\}$ of $G$ by $U_{\hat{\imath}} = \pi^{-1}(\widetilde{U_{\hat{\imath}}})$.  Let $c_{S} \in A$ denote both the barycenter of the face $F_S$ and the exponential of this point in $T \subset G$.  The open set $U_{\hat{\imath}}$ deformation retracts (equivariantly with respect to the conjugation action) to the conjugacy class of $c_{\hat{\imath}}$, which is the quotient $G/Z(c_{\hat{\imath}})$.  More generally, the intersection $\bigcap_{i \notin {S}} U_{\hat{\imath}}$ deformation retracts to the conjugacy class of $c_{S}$, that is to $G/Z(c_S)$. 
We can now describe the $E^1$ term of the twisted Mayer-Vietoris spectral sequence for $K_{\tau}(G)$, based on the cover $\fU$:
\begin{equation} \nn
E^1_{pq} = \bigoplus_{S \subset D, |S|=n-p} K_{\tau, q}(G/Z(c_{S})) \dra K_{\tau,p+q}(G) 
\end{equation}
In an abuse of notation, here the twisting bundle $\tau$ on $G/Z(c_{S})$ is implicitly the restriction of the twisting bundle $\tau$ along the inclusion $\phi: G/Z(c_{S}) \ra G$; this twisting is classified, of course, by the image of the cohomology class $\tau \in H^3(G;\ZZ)$ under the map $\phi^*: H^3(G;\ZZ) \ra H^3(G/Z(c_{S});\ZZ)$.  

This $E^1$ term, and the corresponding $d^1$ differential, can be described more explicitly as follows.  Note that for $S \subset T \subset D$, there is an inclusion $Z(c_S) \subset Z(c_T)$ and a corresponding projection $\pi_{S,T}: G/Z(c_S) \ra G/Z(c_T)$.  The group $G$ is homotopy equivalent to the realization of the simplicial space $Z_k = \bigsqcup_{S \subset D, |S|=n-k} G/Z(c_S)$.  For each $i \in D$, pick a principal $PU(\cH)$ bundle $\tau_i$ on $G/Z(c_{\hat{\imath}})$ and pick isomorphisms $\phi_{ij}: \pi_{\widehat{\imath\jmath},\hat{\imath}}^* \tau_i \xra{\simeq} \pi_{\widehat{\imath\jmath},\hat{\jmath}}^* \tau_j$ such that $\phi_{ij} \phi_{jk} \phi_{ki} = 1$, and such that the resulting principal bundle on $|Z_{\textstyle\cdot}| \simeq G$ is classified by $\tau \in H^3(G;\ZZ)$.  Fix an order $0, 1, \ldots, n$ on the simple affine roots $D$ such that the affine root is first.  For $S \subset D$ let $t(S) \in D$ be the first root that is not in $S$.
The $E^1$ term can now be written
\begin{equation} \nn
E^1_{pq} = \bigoplus_{S \subset D, |S|=n-p} K_{(\pi_{S,\widehat{t(S)}}^* \tau_{t(S)}), q}(G/Z(c_{S})) \dra K_{\tau,p+q}(G) 
\end{equation}
Let $\widehat{S} = \{j_0, j_1, \ldots, j_p\}$ be the complement of $S$, with $j_0 < j_1 < \cdots < j_p$ in the chosen order of the affine roots, and set $T=S \cup j_s$ for some $0 \leq s \leq p$.  With respect to this presentation, the $S$-$T$ component of the $d^1$ differential is $(-1)^s$ times the composite
\begin{equation} \nn
K_{\pi_{S, \widehat{t(S)}}^* \tau_{t(S)}} (G/Z(c_S)) \xra{\phi} K_{\pi_{S,\widehat{t(T)}}^* \tau_{t(T)}} (G/Z(c_S)) = K_{\pi_{S,T}^* \pi_{T, \widehat{t(T)}}^* \tau_{t(T)}} (G/Z(c_S)) \xra{\pi^{S,T}_*} K_{\pi_{T, \widehat{t(T)}}^* \tau_{t(T)}} (G/Z(c_{T})) 
\end{equation}
Here the isomorphism $\phi$ is induced by the chosen twisting isomorphisms $\phi_{ij}$, the map $\pi^{S, T}_*$ is the natural map in twisted $K$-homology, and the homological degree $q$ remains implicit.

This particular presentation of the spectral sequence is unsightly in part because it is meant to isolate the exact role of the twisting isomorphism $\phi$.  In our later computations, we will use that we can work with less standard presentations of the same spectral sequence.

\subsubsection{TMVSS in equivariant $K$-homology} \label{sec-inequivkhom}

We have belabored the non-equivariant twisted Mayer-Vietoris spectral sequence in part because the equivariant spectral sequence is precisely analogous---replacing ``$K$" by ``$K^G$" everywhere in the previous section yields the correct equivariant spectral sequence.  However, establishing that this is the case requires care, because the construction of twisted equivariant $K$-homology is more delicate than that of twisted $K$-homology or twisted equivariant $K$-cohomology.

It is a technical headache to construct a $G$-prespectrum of the homotopy type of the equivariant $K$-theory spectrum that admits an equivariant action by a topological group whose homotopy type is $K(Z,2)$, such that on homotopy the action is given by tensoring an equivariant vector bundle by a line bundle.  Indeed, no such construction exists in the literature, and as a result there is at present no description of a parametrized $G$-prespectrum representing twisted equivariant $K$-homology.  Our primary focus in this note is computational, and as such we do not undertake to build such a parametrized $G$-prespectrum; rather, we sidestep the problem by using a $C^*$-algebra approach to $K$-homology.

Rosenberg~\cite{rosenberg} and Meinrenken~\cite{meinrenken} describe twisted equivariant $K$-homology as follows.  Fix a $G$-Hilbert space $\cH$ containing infinitely many copies of each finite-dimensional $G$-representation.  Equivariant bundles of $G$-$C^*$-algebras over $X$ with fiber the compact operators $\cK_G(\cH)$ and structure group $PU_G(\cH)$ are classified up to Morita equivalence by $H^3_G(X;\ZZ)$.  Such bundles are called ``Dixmier-Douady" bundles.  Pick such a bundle $\cA$ over $X$ whose invariant class is $\tau \in H^3_G(X;\ZZ)$, and for any space $Y$ over $X$, define
$$K^G_{\tau,i}(Y) := KK^i_G(\Gamma(Y,\cA|_Y))$$
where $\Gamma$ is the $G$-$C^*$-algebra of sections vanishing at infinity, and $KK^i_G$ is Kasparov $K$-theory.  In the following, we let $\tau$ refer not only to homology classes in $H^3_G(X;\ZZ)$ but to implicit choices of representing Dixmier-Douady bundles $\cA_\tau$.

This $C^*$-algebra twisted $K$-homology is a generalized homology theory with closed support, that is a generalized Borel-Moore homology, and as such has spectral sequences not for open covers but for closed filtrations.  In particular it has a spectral sequence for the skeletal filtration of the realization $|Z_{\textstyle\cdot}| \simeq G$ of the simplicial space $Z_k = \bigsqcup_{S \subset D, |S|=n-k} G/Z(c_S)$, as follows.

\begin{prop} \label{prop-mvsskh}
For $G$ simple, simply connected, and $\tau \in H^3_G(G;\ZZ)$ represented by a fixed Dixmier-Douady bundle $\cA$, there is a spectral sequence of $K^G(*)$-modules: 
$$E^1_{pq} = \bigoplus_{S \subset D, |S|=n-p} K^G_{\tau, q}(G/Z(c_{S})) \dra K^G_{\tau,p+q}(G)$$
Here $c_{S}$ denotes the barycenter of the face of the Weyl alcove of $\fg$ indexed by the subset $S$ of the affine Dynkin diagram $D$.  

For $S \subset T \subset D$, let $\pi_{S,T}: G/Z(c_S) \ra G/Z(c_T)$ be the projection and let $\widehat{S}$ denote the complement of $S$ in $D$.  Given classes $\tau_i \in H^3_G(G/Z(c_{\hat{\imath}});\ZZ)$ represented by fixed Dixmier-Douady bundles $\cA_i$, together with isomorphisms $\phi_{ij}: \pi_{\widehat{\imath\jmath},\hat{\imath}}^* \cA_i \xra{\simeq} \pi_{\widehat{\imath\jmath},\hat{\jmath}}^* \cA_j$ such that $\phi_{ij} \phi_{jk} \phi_{ki} = 1$, and such that the resulting Dixmier-Douady bundle on $G$ is classified by $\tau \in H^3_G(G;\ZZ)$, the above spectral sequence is presented by
$$E^1_{pq} = \bigoplus_{S \subset D, |S|=n-p} K^G_{(\pi_{S,\widehat{t(S)}}^* \tau_{t(S)}), q}(G/Z(c_{S})) \dra K^G_{\tau,p+q}(G)$$
For $\widehat{S} = \{j_0, j_1, \ldots, j_p\}$ and $T = S \cup j_s$, the $S$-$T$ component of the $d^1$ differential is $(-1)^s$ times the map
$$K^G_{\pi_{S, \widehat{t(S)}}^* \tau_{t(S)}} (G/Z(c_S)) \xra{\phi} K^G_{\pi_{S,\widehat{t(T)}}^* \tau_{t(T)}} (G/Z(c_S)) = K^G_{\pi_{S,T}^* \pi_{T, \widehat{t(T)}}^* \tau_{t(T)}} (G/Z(c_S)) \xra{\pi^{S,T}_*} K^G_{\pi_{T, \widehat{t(T)}}^* \tau_{t(T)}} (G/Z(c_{T}))$$
\end{prop}
\nid This spectral sequence is not new.  It arose originally in the work of Freed, Hopkins, and Teleman, though only a related spectral sequence in $K$-cohomology appears in their papers.  It has been described explicitly by Meinrenken~\cite{meinrenken} and is discussed by Kitchloo and Morava~\cite{kitchloomorava, kitchloo}.  

The main theorem about this particular twisted Mayer-Vietoris spectral sequence in equivariant $K$-homology is the following.
\begin{theorem}[Freed-Hopkins-Teleman~\cite{fht1}] \label{thm-fht}
The twisted Mayer-Vietoris spectral sequence for $K^G_{\tau+h^{\scriptscriptstyle\vee}}(G)$, as in Proposition~\ref{prop-mvsskh}, collapses at the $E^2$ term, and the $E^2$ term, as a $K^G(*)=R[G]$ module, is the fusion ring $F_{\tau}(G)$, that is the ring of positive energy representations of the loop group $LG$ at level $\tau$, concentrated in homological degree $0$.
\end{theorem}
\nid The spectral sequence is 2-periodic in the internal direction, with all groups of odd internal degree vanishing, and this structure is implicit in the statement of the above theorem.  Freed, Hopkins, and Teleman also prove an analogous and much more difficult result for not-necessarily simply connected $G$.

\subsection{A resolution of the fusion ring} \label{sec-resolution}

We redescribe the above twisted Mayer-Vietoris spectral sequence for $K^G_{\tau}(G)$ in terms of twisted representation modules for subgroups of the group $G$, and then even more explicitly in terms of invariants for actions of subgroups of the affine Weyl group of $G$.  In light of Theorem~\ref{thm-fht}, the result is a resolution of $K^G_{\tau}(G)$, therefore of the fusion ring---our resolution differs from that of Meinrenken~\cite{meinrenken} at most by a change of presentation.

The components of the $E^1$ term of our spectral sequence have the form $K^G_{\tau}(G/Z)$, for $Z$ the centralizer of a point of the Weyl alcove.  In order to give a representation-theoretic description of this twisted $K$-homology group, we dualize to twisted $K$-cohomology.  Though $G/Z$ may not be equivariantly $\Spinc$, it is canonically $h^{\scriptscriptstyle\vee}$-twisted equivariantly $\Spinc$, and so satisfies twisted Poincar\'{e} duality:
\begin{equation} \nn
K^G_{\tau+h^{\scriptscriptstyle\vee}}(G/Z) \cong K_G^{\tau}(G/Z) 
\end{equation}
(See~\cite{douglas-tkt} for a definition and discussion of twisted $\Spinc$ structures and~\cite{meinrenken} for a description of the $h^{\scriptscriptstyle\vee}$-twisted structure on adjoint orbits.)  Next we translate the twisted $K$-theory of the quotient $G/Z$ to a twisted $K$-theory of a point: $$K_G^{\tau}(G/Z) \cong K_Z^{\bar{\tau}}(*)$$  
Here $\bar{\tau}$ is the image of $\tau$ under the equivalence $H^3_G(G/Z) = H^3_{Z}(*)$---more precisely, of course, $G$-equivariant Dixmier-Douady bundles on $G/Z$ correspond to $Z$-equivariant Dixmier-Douady bundles on a point, and $\bar{\tau}$ refers to the image under this correspondence of the bundle representing $\tau$.  The group $H^3_{Z}(*)$ classifies $S^1$-central extensions $\widetilde{Z}$ of $Z$, and we can express the twisted $K$-theory in terms of this extension:
$$K_Z^{\bar{\tau},0}(*) = R_{\bar{\tau}}[Z], K_Z^{\bar{\tau},1}(*)=0$$
The ``twisted representation module" $R_{\bar{\tau}}[Z]$ is by definition the submodule of the representation ring of the $\bar{\tau}$-extension $\widetilde{Z}$ consisting of those representations for which the central circle acts by scalar multiplication.  Henceforth we simplify the notation by referring to the $Z$-equivariant Dixmier-Douady bundle on a point corresponding to the $G$-bundle $\tau$ on $G/Z$ not by $\bar{\tau}$ but simply by $\tau$.  

As in Proposition~\ref{prop-mvsskh} we consider any twisting bundle $\tau$ on $G$ as being built from twisting bundles $\tau_i$ on $G/Z(c_{\hat{\imath}})$ via compatible bundle isomorphisms $\phi_{ij}: \pi_{\widehat{\imath\jmath},\hat{\imath}}^* \tau_i \xra{\simeq} \pi_{\widehat{\imath\jmath},\hat{\jmath}}^* \tau_j$.  The map $\pi_{S,T}: G/Z(c_S) \ra G/Z(c_T)$ is the projection corresponding to the inclusion $\iota_{S,T}: Z(c_S) \ra Z(c_T)$ for a pair $S \subset T \subset D$ of subsets of the affine Dynkin diagram $D$.  The above discussion allows us to write the $E^1$ term of the spectral sequence converging to $K^G_{\tau+h^{\scriptscriptstyle\vee}}(G)$ in terms of representation modules:
$$E^1_{p,0} = \bigoplus_{S \subset D, |S|=n-p} R_{\iota_{S,\widehat{t(S)}}^* \tau_{t(S)}}[Z(c_{S})]$$
The spectral sequence is 2-periodic in the internal direction, and henceforth we write only a single homological line.

Next we investigate the $d^1$ differential.  As before let $T = S \cup j_s$ for $j_s \in \widehat{S}$.  In the absence of any twisting, the $S$-$T$ component of the differential in the spectral sequence would be $(-1)^s$ times the natural map
$K^G(G/Z(c_S)) \ra K^G(G/Z(c_T)).$  
Given chosen complex (or more generally $\Spinc$) structures on $G/Z(c_S)$ and $G/Z(c_T)$, therefore chosen Poincar\'{e} duality isomorphisms, this natural map could be described representation-theoretically by holomorphic induction $R[Z(c_S)] \ra R[Z(c_T)]$.  In the presence of twistings, the map is, as one might imagine, a twisted holomorphic induction:
$$R_{\iota_{S, \widehat{t(S)}}^* \tau_{t(S)}} [Z(c_S)] \xra{\phi} R_{\iota_{S,\widehat{t(T)}}^* \tau_{t(T)}} [Z(c_S)] = R_{\iota_{S,T}^* \iota_{T, \widehat{t(T)}}^* \tau_{t(T)}} [Z(c_S)] \xra{\iota^{S,T}_*} R_{\iota_{T, \widehat{t(T)}}^* \tau_{t(T)}} [Z(c_{T})]$$
Here $\iota^{S,T}_*$ is holomorphic induction for the inclusion of central extensions $$\widetilde{Z(c_S)}^{\iota_{S,T}^* \iota_{T, \widehat{t(T)}}^* \tau_{t(T)}} \subset \widetilde{Z(c_T)}^{\iota_{T, \widehat{t(T)}}^* \tau_{t(T)}}$$ with respect to the canonical twisted $\Spinc$ structures on $G/Z(c_S)$ and $G/Z(c_T)$.  The real twist occurs in the map $\phi$, which transforms a representation of $\widetilde{Z(c_S)}^{\iota_{S, \widehat{t(S)}}^* \tau_{t(S)}}$ into a representation of $\widetilde{Z(c_S)}^{\iota_{S,\widehat{t(T)}}^* \tau_{t(T)}}$ according to the bundle isomorphisms $\phi_{ij}$ constructing the twisting $\tau$ from the twistings $\tau_i$.  For example, even when all the twistings $\tau_i$ are trivial, the isomorphisms $\phi_{ij}$ may still involve tensoring by nontrivial line bundles, and the overall differential remains decidedly twisted---this is the type of twisted holomorphic induction that arises in non-equivariant computations of twisted $K$-homology~\cite{douglas-tkt}.  

We can eliminate the central extensions from our description of the spectral sequence by expressing the twisted representation modules as affine Weyl invariants in the representation ring of the maximal torus $T$ of $G$.  For $\alpha \in \ft^*$ a simple root of $G$, let $w_{\alpha}: \ft^* \ra \ft^*$ denote the involution given, as usual, by $w_{\alpha}(\beta) = \beta - \beta(h_{\alpha})\alpha$---this is the reflection in the plane $\ker h_{\alpha}$.  For $\alpha_0 \in \ft^*$ the affine (lowest) root, by slight abuse of notation let $w_{k \alpha_0}: \ft^* \ra \ft^*$ denote the involution given by reflection in the plane $-k \frac{\alpha_0}{2} + \ker h_{\alpha_0}$.  For a subset $S \subset D \backslash \alpha_0$ of the affine Dynkin diagram not containing the affine root, let $W_S^k$ be the group of reflections of $\ft^*$ generated by $\{w_s\}_{s \in S}$---this does not depend on $k$.  For $T = S \cup \alpha_0$, let $W_T^k$ be the group of reflections generated by $W_S^k$ and by $w_{k \alpha_0}$.  The reflection group $W_T^k$ fixes an affine face of the polyhedron bounded by the planes $\{-k \frac{\alpha_0}{2} + \ker h_{\alpha_0}; \ker h_{\alpha_1}; \ker h_{\alpha_2}; \ldots; \ker h_{\alpha_n} \}$; we will refer to this polyhedron as the level $k$ Weyl alcove.

We would like to identify the representation module $R_{\iota_{S,\widehat{t(S)}}^* \tau_{t(S)}}[Z(c_S)]$ in terms of affine Weyl invariants.  For any integer ${k_{t(S)}} \in \ZZ \cong H^3_G(G)$ restricting to the twisting class $\tau_{t(S)} \in H^3_G(G/Z(c_{\widehat{t(S)}}))$, there is an isomorphism $R_{\iota_{S,\widehat{t(S)}}^* \tau_{t(S)}}[Z(c_S)] \cong \ZZ[\Lambda_W]^{W_S^{{k_{t(S)}}}}$---this is seen by comparing the $k_{t(S)}$-affine Weyl invariants in the weight lattice of $G$ to the Weyl invariants in the ${k_{t(S)}}$-affine slice of the weight lattice of the central extension $\widetilde{Z(c_S)}$ associated to the generating twist $1 \in \ZZ \cong H^3_G(G)$.  Let $k \in \ZZ$ correspond to the class of the twisting $\tau \in H^3_G(G)$.  Because the twistings $\tau_i$ glue together to the global twisting $\tau$, we can take ${k_i} = k$ for all $i$.  This choice might seem strange when, for instance, $\tau_i$ is trivial (in which case ${k_i}=0$ would also be a sensible convention), but the differentials in the spectral sequence take a convenient form with respect to this description of the representation modules in terms of $W_S^k$ invariants.  Altogether, the $E^1$ term of the spectral sequence is now
$$E^1_{p,0} = \bigoplus_{S \subset D, |S|=n-p} \ZZ[\Lambda_W]^{W_S^k}$$
This presentation allows us to absorb both the twisting maps $\phi$ and the induction maps $\iota_*$ composing the differential into a single weight induction map, as follows.  For $T = S \cup j_s$ with $j_s \in \widehat{S}$, the $S$-$T$ component of the $d^1$ differential in the spectral sequence is $(-1)^s$ times the map
\begin{align}
\ZZ[\Lambda_W]^{W_S^k} &\ra \ZZ[\Lambda_W]^{W_T^k} \nn\\
\left[\frac{A_{\mu + \rho_S}^{W_S^k}}{A_{\rho_S}^{W_S^0}}\right] &\mapsto
\left[\frac{A_{\mu + \rho_T}^{W_T^k}}{A_{\rho_T}^{W_T^0}}\right] \nn
\end{align}
Here $\mu$ is a weight in $\Lambda_W$, and $\rho_S$ is the half sum of the positive roots of $Z(c_S)$, with respect to the Weyl chamber determined by the subset $S$ of the affine Dynkin diagram $D$---similarly for $\rho_T$.  For a weight $\lambda \in \Lambda_W$, the expression $A_{\lambda}^W$ denotes the antisymmetrization of $\lambda$ with respect to the reflection group $W$.  Note that the invariant module $\ZZ[\Lambda_W]^{W_S^k}$ is spanned by the quotients $\frac{A_{\mu + \rho_S}^{W_S^k}}{A_{\rho_S}^{W_S^0}}$.  Of course, this induction map is just the usual weight description of holomorphic induction, adjusted by an affine translation depending on the twisting $k$.
Using Theorem~\ref{thm-fht} we can summarize this presentation of the twisted Mayer-Vietoris spectral sequence in twisted equivariant $K$-homology as follows.
\begin{prop} \label{prop-resolution}
The complex of $R[G]$-modules
$\bigoplus_{S \subset D, |S|=n-p} \ZZ[\Lambda_W]^{W_S^k}$, with differential having $S$-$T$ component $\left[\frac{A_{\mu + \rho_S}^{W_S^k}}{A_{\rho_S}^{W_S^0}}\right] \mapsto (-1)^s \left[\frac{A_{\mu + \rho_T}^{W_T^k}}{A_{\rho_T}^{W_T^0}}\right]$, where $T=S \cup j_s$, is acyclic except in degree zero, where it has homology the level $k$ fusion ring $F_k(G)$.
\end{prop}

\section{The fusion ring of $G_2$} \label{sec-g2}

The complex coming from the Mayer-Vietoris spectral sequence in twisted equivariant $K$-homology provides a presentation of the fusion ring.  We illustrate such presentations by explicitly describing the fusion ring of $G_2$.

\begin{theorem} \label{thm-g2}
For $k>0$, the level $k$ fusion ring of $G_2$ is given as an $R[G_2]$-module, and therefore as a ring, by
\begin{equation} \nn
F_k[G_2] =
\begin{cases}
R[G_2]/(\rho_{(1, \frac{k}{2})}; \rho_{(0, \frac{k}{2})}+\rho_{(0,\frac{k}{2} + 1)}; \rho_{(1,\frac{k}{2} - 1)}+\rho_{(1,\frac{k}{2} +1)}; \rho_{(k+2,0)})
& \text{$k$ even} \\
R[G_2]/(\rho_{(0, \frac{k+1}{2})};\rho_{(1, \frac{k-1}{2})}+\rho_{(1,\frac{k+1}{2})}; \rho_{(0, \frac{k-1}{2})}+\rho_{(0,\frac{k+3}{2})}; \rho_{(k+2,0)})
& \text{$k$ odd}
\end{cases}
\end{equation}
Here $\rho_{(a,b)}$ is the irreducible representation of $G_2$ with highest weight $\lambda_{1,0}^a \lambda_{0,1}^b$, for $\lambda_{1,0}$ and $\lambda_{0,1}$ respectively the short and long fundamental weights of $G_2$.
\end{theorem}

\begin{proof}

The $E^1$ term of the twisted Mayer-Vietoris spectral sequence has the form
$$ R[G_2] \oplus R_k[SO(4)] \oplus R[SU(3)] \xla{d^{\tw}} R[U(2)_v] \oplus R[U(2)_d] \oplus R[U(2)_h] \xla{d^{\tw}} R[T] $$
Here $d^{\tw}$ is twisted holomorphic induction.  Further $U(2)_v$, $U(2)_d$, and $U(2)_h$ denote the subgroups of $G_2$ whose representations rings are $\ZZ[\Lambda_W]^{W_v}$, $\ZZ[\Lambda_W]^{W_d}$, and $\ZZ[\Lambda_W]^{W_h}$, where $W_v$, $W_d$, and $W_h$ are the reflections in the vertical, diagonal, and horizontal walls of the Weyl chamber decomposition of the dual of the Lie algebra of the torus of $G_2$.  Note that $U(2)_d$ and $U(2)_h$ are generated by pairs of long roots of $G_2$, while $U(2)_v$ is generated by a pair of short roots.  Also, up to $R[G_2]$-module isomorphism the twisted representation module $R_k[SO(4)]$ only depends on the parity of the level $k$.

We first compute the various summands of this $E^1$ term as $R[G_2]$-modules.
\begin{lemma} \label{lemma-basesg2}
The representation rings $R[T]$, $R[U(2)_s]$, $R[U(2)_l]$, $R[SU(3)]$, and $R[SO(4)]$, and the twisted representation module $R_1[SO(4)]$ have the following structure, as $R:=R[G_2]$-modules.  Here $U(2)_s$, respectively $U(2)_l$, is the $U(2)$ subgroup of $G_2$ generated by the eigenspaces for $\pm \alpha$ for $\alpha$ a short, respectively long, simple root of $G_2$.
\begin{align*}
R[T] &= \: \scriptstyle{R \, \rho^T_{0,0} \oplus R \, \rho^T_{0,1} \oplus R \, \rho^T_{0,2} \oplus R \, \rho^T_{1,0} \oplus R \, \rho^T_{1,1} \oplus R \, \rho^T_{1,2} \oplus R \, \rho^T_{2,-1} \oplus R \, \rho^T_{2,0} \oplus R \, \rho^T_{2,1} \oplus R \, \rho^T_{3,-1} \oplus R \, \rho^T_{3,0} \oplus R \, \rho^T_{3,1}} \\
R[U(2)_s] &= R \, \rho^{U(2)_s}_{0,0} \oplus R \, \rho^{U(2)_s}_{0,1} \oplus R \, \rho^{U(2)_s}_{0,2} \oplus R \, \rho^{U(2)_s}_{1,0} \oplus R \, \rho^{U(2)_s}_{1,1} \oplus R \, \rho^{U(2)_s}_{1,2} \\
R[U(2)_l] &= R \, \rho^{U(2)_l}_{0,0} \oplus R \, \rho^{U(2)_l}_{1,0} \oplus R \, \rho^{U(2)_l}_{2,0} \oplus R \, \rho^{U(2)_l}_{3,0} \oplus R \, \rho^{U(2)_l}_{0,1} \oplus R \, \rho^{U(2)_l}_{-2,2} \\
R[SU(3)] &= R \, \rho^{SU(3)}_{0,0} \oplus R \, \rho^{SU(3)}_{-1,0} \\
R[SO(4)] &= R \, \rho^{SO(4)}_{0,0} \oplus R \, \rho^{SO(4)}_{1,-1} \oplus R \, \rho^{SO(4)}_{0,-1} \\
R_1[SO(4)] &= R \, \rho^{SO(4),1}_{0,0} \oplus R \, \rho^{SO(4),1}_{1,0} \oplus R \, \rho^{SO(4),1}_{1,-1}
\end{align*}
The notation $\rho^H_{a,b}$ refers to the irreducible representation of $H$ with highest weight $\lambda_{1,0}^a \lambda_{0,1}^b$; as before, $\lambda_{1,0}$ and $\lambda_{0,1}$ are respectively the short and long fundamental weights of $G_2$.  The element $\rho^{SO(4),1}_{a,b} \in \ZZ[\Lambda_W]^{W_{SO(4)}^1}$ is the irreducible representation of the generating central extension of $SO(4)$ with the indicated highest weight.  The Weyl chambers for all groups have been determined by the corresponding subsets of the affine simple roots.
\end{lemma} 

\begin{proof}
The first five presentations are determined by explicit, and tedious, computation.  For example, note that $R[SO(4)] \cong \ZZ[m,n,p]/(p^2-mn-m-n-1)$ with $R[G_2] \cong \ZZ[a,b]$ action by $a = m+p$ and $b=pm-p+n+m$; the elements $\{1,p,n\}$ provide a module basis; here $m$, $n$, and $p$ are respectively the irreducible representations of $SO(4)$ with highest weights $(2,-1)$, $(0,-1)$, and $(1,-1)$.  (It is interesting to note that though $R[SO(4)]$ is free over $R[G_2]$ it is itself a singular variety.)  Note that the basis we have chosen for $R[U(2)_l]$ is carefully tailored for use in the proof of Theorem~\ref{thm-g2}, and does not arise in a natural way---for instance, we do not know if it is a Gr\"{o}bner basis for any presentation of the module.

The calculation for the twisted representation module $R_1[SO(4)]$ relies on the untwisted case.  As an $R[SO(4)]$-module $R_1[SO(4)]=\ZZ[\Lambda_W]^{W_{SO(4)}^1}$ is generated by the two 2-dimensional irreducible modules $r$ and $s$, whose highest weights are $(1,0)$ and $(0,0)$ respectively.  In light of the presentation of $R[SO(4)]$, this implies that $R_1[SO(4)]$ is generated over $R[G_2]$ by $\{r,s,r^2 s, s^2 r, s^3\}$.  Note that the generating representations of $G_2$ are $a=rs+r^2-1$ and $b=r^3 s - 2rs + r^2 +s^2 - 2$, and that $as=s^2 r + r^2 s - s$ and $-a(r^2 s - 2s+r)+b(r+s)=r^2 s - 4s - r + rs^2 + s^3$.  The indicated basis follows.
\end{proof}  

Proposition~\ref{prop-resolution} presents the $E^1$ term of the twisted Mayer-Vietoris spectral sequence as a complex whose homology is the fusion ring, with terms given by affine Weyl invariants:
$$\ZZ[\Lambda_W]^{W_{G_2}} \oplus \ZZ[\Lambda_W]^{W_{SO(4)}^k} \oplus \ZZ[\Lambda_W]^{W_{SU(3)}^k} \xla{d_1} \ZZ[\Lambda_W]^{W_{U(2)_v}} \oplus \ZZ[\Lambda_W]^{W_{U(2)_d}} \oplus \ZZ[\Lambda_W]^{W_{U(2)_h}^k} \xla{d_2} \ZZ[\Lambda_W]^{W_{T}}$$
Substituting the bases from Lemma~\ref{lemma-basesg2} into this presentation, the complex in question appears as in Figure~\ref{fig-g2}.
\begin{center}
\begin{figure}[ht]
\centering\includegraphics[height=3.5in]{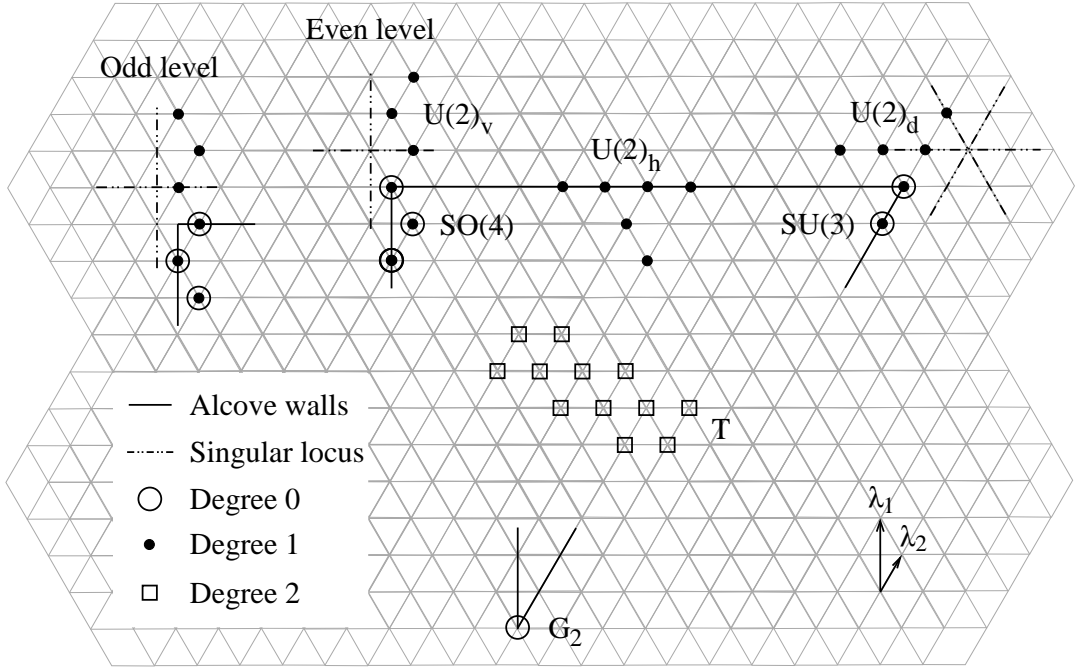}
\caption{Basis of highest weights for an $R[G_2]$-module resolution of the fusion ring $F_k[G_2]$.} \label{fig-g2}
\end{figure}
\end{center}

\vspace*{-10pt}
To compute the desired cokernel of $d_1$, observe that the induction map $d_2: \ZZ[\Lambda_W]^{W_{T}} \xra{d_2} \ZZ[\Lambda_W]^{W_{U(2)_h}^k}$ is surjective.  This implies that $$d_1(\ZZ[\Lambda_W]^{W_{U(2)_h}^k}) \subset d_1(\ZZ[\Lambda_W]^{W_{U(2)_v}} \oplus \ZZ[\Lambda_W]^{W_{U(2)_d}})$$
Next consider $d_1(\ZZ[\Lambda_W]^{W_{U(2)_d}})$.  The representations $\rho^{U(2)_d}_{(k,0)}$ and $\rho^{U(2)_d}_{(k-1,0)}$ in $\ZZ[\Lambda_W]^{W_{U(2)_d}}$ induce respectively to the irreducible generators $\rho^{SU(3),k}_{(k,0)}$ and $\rho^{SU(3),k}_{(k-1,0)}$ of $\ZZ[\Lambda_W]^{W_{SU(3)}^k}$.  A priori, the remaining four generators of $\ZZ[\Lambda_W]^{W_{U(2)_d}}$ would create four relations in the fusion ideal.  However, three of the generators of the diagonal Weyl invariants appear on the singular wall for holomorphic induction from the representation ring of the torus to the affine representation ring of the horizontal $U(2)$.  As a result, the differential has the form
\begin{align}
d_2(\rho^T_{(k+1,0)}) &= (V, -\rho^{U(2)_d}_{(k+1,0)}, 0) \nn\\
d_2(\rho^T_{(k-1,1)}) &= (V', -\rho^{U(2)_d}_{(k-1,1)}, 0) \nn\\
d_2(\rho^T_{(k-3,2)}) &= (V'', -\rho^{U(2)_d}_{(k-3,2)}, 0) \nn
\end{align}
where $V$, $V'$, and $V''$ are elements of $\ZZ[\Lambda_W]^{W_{U(2)_v}}$.  It follows that
$$d_1(\langle \rho^{U(2)_d}_{(k+1,0)}, \rho^{U(2)_d}_{(k-1,1)}, \rho^{U(2)_d}_{(k-3,2)} \rangle) \subset d_1(\ZZ[\Lambda_W]^{W_{U(2)_v}})$$  The remaining generator of $\ZZ[\Lambda_W]^{W_{U(2)_d}}$ has differential
$$d_1(\rho^{U(2)_d}_{(k+2,0)}) = (-\rho^{G_2}_{(k+2,0)},0) \in \ZZ[\Lambda_W]^{W_{G_2}} \oplus \ZZ[\Lambda_W]^{W_{SU(3)}^k}$$
Finally, consider $d_1(\ZZ[\Lambda_W]^{W_{U(2)_v}})$.  If the level $k$ is even, the representations $\rho^{U(2)_v}_{(0,\frac{k}{2})}$, $\rho^{U(2)_v}_{(1,\frac{k}{2}-1)}$, and $\rho^{U(2)_v}_{(0,\frac{k}{2}-1)}$ induce on the one hand to the generators $\rho^{SO(4),k}_{(0,\frac{k}{2})}$, $\rho^{SO(4),k}_{(1,\frac{k}{2}-1)}$, and $\rho^{SO(4),k}_{(0,\frac{k}{2}-1)}$ of $\ZZ[\Lambda_W]^{W_{SO(4)}^k}$ and on the other hand to the corresponding irreducible representations in $\ZZ[\Lambda_W]^{W_{G_2}}$.  The remaining generators of $\ZZ[\Lambda_W]^{W_{U(2)_v}}$ have differentials
\begin{align}
d_1(\rho^{U(2)_v}_{(1,\frac{k}{2})}) &= (-\rho^{G_2}_{(1, \frac{k}{2})}, 0) \nn\\
d_1(\rho^{U(2)_v}_{(0,\frac{k}{2}+1)}) &= (-\rho^{G_2}_{(0,\frac{k}{2} + 1)}, -\rho^{SO(4),k}_{(0,\frac{k}{2})}) \nn\\
d_1(\rho^{U(2)_v}_{(1,\frac{k}{2}+1)}) &= (-\rho^{G_2}_{(1,\frac{k}{2} +1)}, -\rho^{SO(4),k}_{(1,\frac{k}{2}-1)}) \nn
\end{align}
If the level $k$ is odd, the representations $\rho^{U(2)_v}_{(1,\frac{k-1}{2})}$, $\rho^{U(2)_v}_{(0,\frac{k-1}{2})}$, and $\rho^{U(2)_v}_{(1,\frac{k-3}{2})}$ induce to the corresponding generators of $\ZZ[\Lambda_W]^{W_{SO(4)}^k}$ and to the corresponding irreducible representations in $\ZZ[\Lambda_W]^{W_{G_2}}$.  The remaining generators of $\ZZ[\Lambda_W]^{W_{U(2)_v}}$ have differentials
\begin{align}
d_1(\rho^{U(2)_v}_{(0,\frac{k+1}{2})}) &= (-\rho^{G_2}_{(0,\frac{k+1}{2})}, 0) \nn\\
d_1(\rho^{U(2)_v}_{(1,\frac{k+1}{2})}) &= (-\rho^{G_2}_{(1,\frac{k+1}{2})}, -\rho^{SO(4),k}_{(1,\frac{k-1}{2})}) \nn\\
d_1(\rho^{U(2)_v}_{(0,\frac{k+3}{2})}) &= (-\rho^{G_2}_{(0,\frac{k+3}{2})}, -\rho^{SO(4),k}_{(0,\frac{k-1}{2})}) \nn
\end{align}
The indicated presentations of the fusion ring follow.
\end{proof}

The fusion ideal of $G$ at level $k$, that is the ideal $I_k \subset R[G]$ such that $F_k[G] = R[G]/I_k$, contains and is often generated by the irreducible representations of $G$ whose highest weights are on the affine wall of the level $k+1$ Weyl alcove; the order of this collection of irreducible representations tends to infinity with the level.  However, the above theorem shows that for $G_2$, the fusion ideal is generated by 4 representations, independent of the level---in section~\ref{sec-finiteness} we will see that such a finite bound exists for all groups.  Note that we expect that a minimal presentation of the $G_2$ fusion ring can be obtained from the presentation in Theorem~\ref{thm-g2}, at both even and odd levels, by simply omitting the relation $\rho_{(k+2,0)}$.

For a group $G$ of rank $n$, the fusion ideal $I_k$ has at least $n$ generators.  When this lower bound is achieved, the fusion ring is a zero-dimensional complete intersection.  The complexification of the fusion ring $F_k[G] \otimes \CC$ always has this form, but in general the structure of the integral fusion ring is more subtle.

\section{Finiteness of the fusion ideal} \label{sec-finiteness}

The fusion ring of positive energy representation of the loop group $LG$ is the homology of a complex whose terms are sums of submodules of representation rings of central extensions of subgroups of $G$.  In this section we describe the structure of these so-called twisted representation modules, proving in particular that they are always free as modules over the representation ring of $G$.  We then generalize the $G_2$ computation from section~\ref{sec-g2} to all simple simply connected groups, proving that there is a level-independent bound on the number of generators of the fusion ideal.

\subsection{Twisted representation modules} \label{sec-twrepmod}

In section~\ref{sec-resolution}, we described a resolution of the fusion ring $F_k[G]$ whose terms had the form of twisted representation modules $R_k[H]$.  Recall that the group $G$ is simply connected, the subgroup $H$ is the centralizer in $G$ of an element of the Weyl alcove of $G$, the central extension $\widetilde{H}$ is classified by the image of the generating twist of $G$ under the map $H^3_G(G) \ra H^3_G(G/H) = H^3_H(*)$, and the $R[H]$-module $R_k[H]$ is the submodule of representations of $\widetilde{H}$ for which the central circle acts by the $k$-th power of scalar multiplication.  In this section we describe $R_k[H]$ as an $R[G]$-module, and enumerate all the subgroups $H$ of all groups $G$ for which there exists a $k$ for which the module $R_k[H]$ is indeed twisted, that is not isomorphic to the representation ring $R[H]$.  All semi-simple subgroups of $E_8$ arising as centralizers are twisted in this sense, and for entertainment we describe the structure of these groups in detail.

Modulo the Serre conjecture (that is the Quillen-Suslin theorem), Pittie~\cite{pittie} proved that for $G$ a compact connected simply connected Lie group and $H$ a closed connected subgroup of maximal rank, the representation ring $R[H]$ is a free $R[G]$-module.  The same techniques are effective for studying twisted representation modules:

\begin{prop} \label{prop-free}
Let $G$ be compact, connected, simply connected, simple, and let $H$ be the centralizer of an element of $G$.  For all $k$ the twisted representation module $R_k[H]$ is a free $R[G]$-module.
\end{prop}
\begin{proof} 
The central extension $\widetilde{H}$ is the circle bundle associated to a connected principal cyclic-group subbundle $\bar{H} \subset \widetilde{H}$.  Let $\bar{T}$ denote the maximal torus of $\bar{H}$.  The twisted representation module $R_k[H]$ is an $R[H]$-submodule of $R[\bar{H}]$, and similarly $R_k[T]$ is an $R[H]$-submodule of $R[\bar{T}]$.  Here the modules $R_k[H]$ and $R_k[T]$ are the submodules of level $k$ representations of the extensions $\widetilde{H}$ and $\widetilde{T}$; these extensions are defined by the restriction of the generating twist of $G$ to the orbits $G/H$ and $G/T$.  The inclusions $R_k[H] \ra R[\bar{H}]$ and $R_k[T] \ra R[\bar{T}]$ commute with both the restriction and the induction maps $R[\bar{H}] \leftrightarrows R[\bar{T}]$ and $R_k[H] \leftrightarrows R_k[T]$.  The induction map from $R[\bar{T}]$ to $R[\bar{H}]$ splits the corresponding restriction, and as a result the induction map from $R_k[T]$ to $R_k[H]$ also splits the restriction.  The module $R_k[H]$ is thereby an $R[H]$-direct summand of $R_k[T]$.  The module $R_k[T]$ is isomorphic as an $R[H]$-module to the untwisted module $R[T]$, which in turn is free as an $R[G]$-module.  Hence the module $R_k[H]$ is a projective, therefore by the Quillen-Suslin theorem free $R[G]$-module, as desired.
\end{proof}

\nid No doubt, twisted representation modules are free in greater generality than this proposition describes.  Note that the module $R_k[H]$ is free of the same rank as $R[H]$, namely $|W_G|/|W_H|$.

The splitting described in the above proof can also be seen in terms of affine Weyl invariants.  The representation ring $R[H]$ is isomorphic to the Weyl invariants $\ZZ[\Lambda_W]^{W_H}$ in the weight lattice, and the twisted representation module $R_k[H]$ is isomorphic, as an $R[H]$-module, to the invariants $\ZZ[\Lambda_W]^{W_H^k}$; here the affine Weyl action ``$W_H^k$" is defined, as in section~\ref{sec-resolution}, by reflections in the Weyl hyperplanes containing a particular face of the level $k$ Weyl alcove.  Let $\bar{\Lambda}_W \subset \ft^*$ denote the weight lattice of the covering group $\bar{H}$ of $H$---as a lattice in $\ft^*$, it inherits the various Weyl actions.  There is a commutative diagram
\begin{equation} \nn
\xymatrix{
\ZZ[\Lambda_W]^{W_H^k} \ar@{->}[r] \ar@<-.5ex>@{->}[d] & \ZZ[\bar{\Lambda}_W]^{W_H^k} \ar[r]^{\cdot \mu^{-1}} \ar@<-.5ex>@{->}[d] &
\ZZ[\bar{\Lambda}_W]^{W_H^0} \ar@<-.5ex>@{->}[d] \\
\ZZ[\Lambda_W] \ar@<-.5ex>@{-->}[u] \ar@{->}[r] & \ZZ[\bar{\Lambda}_W] \ar[r]^{\cdot \mu^{-1}} \ar@<-.5ex>@{-->}[u] &
\ZZ[\bar{\Lambda}_W] \ar@<-.5ex>@{-->}[u]
}
\end{equation}
The upward maps are affine induction to the invariant modules.  The righthand horizontal maps are multiplication by $\mu^{-1}$ for some weight $\mu \in (\bar{\Lambda}_W)^{W_H^k}$ fixed by the affine action.  The rightmost induction map splits the corresponding inclusion, and so the leftmost induction must similarly provide a splitting.

We now briefly describe when $R_k[H]$ is indeed twisted, for the various centralizer subgroups of the simple simply connected groups $G$.  This information is essential for explicit computations of the fusion rings, as illustrated for $G_2$ in section~\ref{sec-g2}.  As before, for a subset $S$ of the affine Dynkin diagram $D$ of $G$, the corresponding group $H_S$ is the centralizer of the face $F_S$ of the Weyl alcove, where $F_S$ is defined as the fixed set of the subgroup of the affine Weyl group generated by the reflections corresponding to the simple roots $\{\alpha_i\}_{i \in S}$.

\begin{prop} 
For a subset $S$ of the affine Dynkin diagram $D$ of $G$, the representation module $R_k[H_S]$ is isomorphic to $R[H_S]$ as an $R[H_S]$-module precisely when $\gcd_{i \in \widehat{S}} h_i^{\scriptscriptstyle\vee}$ divides the level $k$.
\end{prop}  
\begin{proof}
Recall that the Coxeter labels $h_i \in \NN$ of the nodes $\alpha_i$ of the nonaffine Dynkin diagram $\bar{D}$ are defined by $-\alpha_0 = \sum_{i=1}^n h_i \alpha_i$.  The dual Coxeter labels $h_i^{\scriptscriptstyle\vee}$ are defined in terms of the Coxeter labels by $h_i^{\scriptscriptstyle\vee} = h_i F(\alpha_i, \alpha_i) / F(\alpha_0, \alpha_0)$, where $F$ is an invariant inner product.  Let $\{b_i\} \subset \ft$ be the coroots, that is the unique vectors in the tori of the fundamental algebras $\mathfrak{su}(2)_i$ such that $\alpha_i(b_i)=2$.  The fundamental weights $\lambda_i$ are by definition evaluation dual to the coroots.  Identify $\ft$ and $\ft^*$ using the inner product $F$.

Observe that the plane spanned by the face $F_S^k$ of the level $k$ Weyl alcove contains a weight of the torus $T$ of $G$ precisely when the twisted representation module $R_k[H_S]$ is isomorphic to $R[H_S]$ as an $R[H_S]$-module.  When $S=\hat{\imath} \subset D$ is a subset of the affine Dynkin diagram of order $n$, the vertex $F_S^k$ is a weight of $T$ when $t_i := -2 F(\alpha_0, \lambda_i) / F(\alpha_0, \alpha_0)$ divides the level $k$.  Note that
$$t_i = \frac{-2 F(\alpha_0, \lambda_i)}{F(\alpha_0, \alpha_0)} = \frac{2 F(\sum_{j=1}^n h_j \alpha_j, \lambda_i)}{F(\alpha_0, \alpha_0)} = \sum_{j=1}^n \frac{h_j F(b_j, \lambda_i) F(\alpha_i, \alpha_i)}{F(\alpha_0, \alpha_0)} = \frac{h_i F(\alpha_i, \alpha_i)}{F(\alpha_0, \alpha_0)} = h_i^{\scriptscriptstyle\vee}$$
More generally, for any subset $S \subset D$, the span of the face $F_S^k$ contains a weight when the greatest common divisor $\gcd_{i \in \widehat{S}} h_i^{\scriptscriptstyle\vee}$ divides the level.
\end{proof}
\pagebreak

All the dual Coxeter labels for groups of type A and type C are 1.  There are no twisted representation modules for these groups, and the resulting fusion rings are simpler as a result.  The dual Coxeter labelings in types B and D are

\hspace*{-30pt}
\begin{minipage}[c]{.5\textwidth}
\begin{center}
\begin{equation} \nn
\xymatrix@-10pt{ 
& \scriptstyle \alpha_0 \ar@{-}[d] \save[]+<-11pt,0pt>*{\scriptstyle 1} \restore \\ \scriptstyle
\alpha_1 \ar@{-}[r] \save[]+<0pt,-10pt>*{\scriptstyle 1} \restore & \scriptstyle \alpha_2 \ar@{-}[r] \save[]+<0pt,-10pt>*{\scriptstyle 2} \restore & \scriptstyle \alpha_3 \ar@{-}[r] \save[]+<0pt,-10pt>*{\scriptstyle 2} \restore & \cdots \ar@{-}[r] & \scriptstyle \alpha_{n-1} \ar@{=}[r] |-{\SelectTips{cm}{12}\object@{>}} \save[]+<0pt,-10pt>*{\scriptstyle 2} \restore & \scriptstyle \alpha_n \save[]+<0pt,-10pt>*{\scriptstyle 1} \restore 
}
\end{equation}
\end{center}
\end{minipage}
\begin{minipage}[c]{.5\textwidth}
\begin{center}
\begin{equation} \nn
\xymatrix@-10pt{
& \scriptstyle \alpha_0 \ar@{-}[d] \save[]+<-11pt,0pt>*{\scriptstyle 1} \restore &&&& \scriptstyle \alpha_n \ar@{-}[d] \save[]+<-11pt,0pt>*{\scriptstyle 1} \restore \\
\scriptstyle \alpha_1 \ar@{-}[r] \save[]+<0pt,-10pt>*{\scriptstyle 1} \restore & \scriptstyle \alpha_2 \ar@{-}[r] \save[]+<0pt,-10pt>*{\scriptstyle 2} \restore & \scriptstyle \alpha_3 \ar@{-}[r] \save[]+<0pt,-10pt>*{\scriptstyle 2} \restore & \cdots \ar@{-}[r] & \scriptstyle \alpha_{n-3} \ar@{-}[r] \save[]+<0pt,-10pt>*{\scriptstyle 2} \restore & \scriptstyle \alpha_{n-2} \ar@{-}[r] \save[]+<0pt,-10pt>*{\scriptstyle 2} \restore & \scriptstyle \alpha_{n-1} \save[]+<0pt,-10pt>*{\scriptstyle 1} \restore 
}
\end{equation}
\end{center}
\end{minipage}
\vspace{8pt}

\nid In type B, the module $R_k[H_S]$ is twisted when $k$ is odd and $S$ contains $\alpha_0$, $\alpha_1$, and $\alpha_n$.  
In type D, the module $R_k[H_S]$ is twisted when $k$ is odd and $S$ contains $\alpha_0$, $\alpha_1$, $\alpha_{n-1}$, and $\alpha_n$.

The exceptional dual Coxeter labels are

\hspace*{-30pt}
\begin{minipage}[c]{.5\textwidth}
\begin{center}
\begin{equation} \nn
\xymatrix@-10pt{
\scriptstyle \alpha_0 \ar@{-}[r] \save[]+<0pt,-10pt>*{\scriptstyle 1} \restore & \scriptstyle \alpha_1 \ar@3{-}[r] |-{\SelectTips{cm}{12}\object@{>}} \save[]+<0pt,-10pt>*{\scriptstyle 2} \restore 
& \scriptstyle \alpha_2 \save[]+<0pt,-10pt>*{\scriptstyle 1} \restore
}
\end{equation}
\end{center}
\end{minipage}
%
%
\begin{minipage}[c]{.5\textwidth}
\begin{center}
\begin{equation} \nn
\xymatrix@-10pt{
\scriptstyle \alpha_0 \ar@{-}[r] \save[]+<0pt,-10pt>*{\scriptstyle 1} \restore & \scriptstyle \alpha_1 \ar@{-}[r] \save[]+<0pt,-10pt>*{\scriptstyle 2} \restore & \scriptstyle \alpha_2 \ar@{=}[r] |-{\SelectTips{cm}{12}\object@{>}} \save[]+<0pt,-10pt>*{\scriptstyle 3} \restore 
& \scriptstyle \alpha_3 \save[]+<0pt,-10pt>*{\scriptstyle 2} \restore \ar@{-}[r] & \scriptstyle \alpha_4 \save[]+<0pt,-10pt>*{\scriptstyle 1} \restore
}
\end{equation}
\end{center}
\end{minipage}

\vspace{5pt}
\hspace*{-30pt}
\begin{minipage}[c]{.5\textwidth}
\begin{center}
\begin{equation} \nn
\xymatrix@-10pt{
&& \scriptstyle \alpha_0 \ar@{-}[d] \save[]+<-11pt,0pt>*{\scriptstyle 1} \restore && \\
&& \scriptstyle \alpha_2 \ar@{-}[d] \save[]+<-11pt,0pt>*{\scriptstyle 2} \restore && \\
\scriptstyle \alpha_1 \ar@{-}[r] \save[]+<0pt,-10pt>*{\scriptstyle 1} \restore &
\scriptstyle \alpha_3 \ar@{-}[r] \save[]+<0pt,-10pt>*{\scriptstyle 2} \restore &
\scriptstyle \alpha_4 \ar@{-}[r] \save[]+<0pt,-10pt>*{\scriptstyle 3} \restore &
\scriptstyle \alpha_5 \ar@{-}[r] \save[]+<0pt,-10pt>*{\scriptstyle 2} \restore &
\scriptstyle \alpha_6 \save[]+<0pt,-10pt>*{\scriptstyle 1} \restore
}
\end{equation}
\end{center}
\end{minipage}
\begin{minipage}[c]{.5\textwidth}
\begin{center}
\begin{equation} \nn
\xymatrix@-10pt{
&&& \scriptstyle \alpha_2 \ar@{-}[d] \save[]+<-11pt,0pt>*{\scriptstyle 2} \restore \\
\scriptstyle \alpha_0 \ar@{-}[r] \save[]+<0pt,-10pt>*{\scriptstyle 1} \restore &
\scriptstyle \alpha_1 \ar@{-}[r] \save[]+<0pt,-10pt>*{\scriptstyle 2} \restore &
\scriptstyle \alpha_3 \ar@{-}[r] \save[]+<0pt,-10pt>*{\scriptstyle 3} \restore &
\scriptstyle \alpha_4 \ar@{-}[r] \save[]+<0pt,-10pt>*{\scriptstyle 4} \restore &
\scriptstyle \alpha_5 \ar@{-}[r] \save[]+<0pt,-10pt>*{\scriptstyle 3} \restore &
\scriptstyle \alpha_6 \ar@{-}[r] \save[]+<0pt,-10pt>*{\scriptstyle 2} \restore &
\scriptstyle \alpha_7 \save[]+<0pt,-10pt>*{\scriptstyle 1} \restore
}
\end{equation}

\vspace{0pt}

\begin{equation} \nn
\xymatrix@-10pt{
&& \scriptstyle \alpha_2 \ar@{-}[d] \save[]+<-11pt,0pt>*{\scriptstyle 3} \restore \\
\scriptstyle \alpha_1 \ar@{-}[r] \save[]+<0pt,-10pt>*{\scriptstyle 2} \restore &
\scriptstyle \alpha_3 \ar@{-}[r] \save[]+<0pt,-10pt>*{\scriptstyle 4} \restore &
\scriptstyle \alpha_4 \ar@{-}[r] \save[]+<0pt,-10pt>*{\scriptstyle 6} \restore &
\scriptstyle \alpha_5 \ar@{-}[r] \save[]+<0pt,-10pt>*{\scriptstyle 5} \restore &
\scriptstyle \alpha_6 \ar@{-}[r] \save[]+<0pt,-10pt>*{\scriptstyle 4} \restore &
\scriptstyle \alpha_7 \ar@{-}[r] \save[]+<0pt,-10pt>*{\scriptstyle 3} \restore &
\scriptstyle \alpha_8 \ar@{-}[r] \save[]+<0pt,-10pt>*{\scriptstyle 2} \restore &
\scriptstyle \alpha_0 \save[]+<0pt,-10pt>*{\scriptstyle 1} \restore
}
\end{equation}
\end{center}
\end{minipage}
\vspace{8pt}
 
\nid Let $\langle\fh\rangle$ refer to the subgroup of $G$ whose Lie algebra is generated by $\fh$ and $\ft$.  The only twisted representation module for $G_2$ is, as we saw in section~\ref{sec-g2}, the module $R_k[\langle\mathfrak{su}(2) \times \mathfrak{su}(2)\rangle]$ for odd level $k$.  For the group $F_4$, the modules $R_k[\langle\mathfrak{sp}(3) \times \mathfrak{su}(2)\rangle]$ and $R_k[\langle\mathfrak{su}(4) \times \mathfrak{su}(2)\rangle]$ are twisted for odd $k$, while $R_k[\langle\mathfrak{su}(3)^2\rangle]$ is twisted for $k$ not divisible by 3; the only other twisted module is $R_k[\langle\mathfrak{su}(2)^3\rangle]$ at odd level.  For $E_6$, the three modules of the form $R_k[\langle\mathfrak{su}(6) \times \mathfrak{su}(2)\rangle]$ are twisted at odd level, while $R_k[\langle\mathfrak{su}(3)^3\rangle]$ is twisted when $k$ is not divisible by 3; similarly, three modules of the form $R_k[\langle\mathfrak{su}(4) \times \mathfrak{su}(2)^2\rangle]$ and the module $R_k[\langle\mathfrak{su}(2)^4\rangle]$ are twisted at odd level.  For $E_7$, we encounter the first order 4 twisting, for $R_k[\langle\mathfrak{su}(4)^2 \times \mathfrak{su}(2)\rangle]$; the two modules $R_k[\langle\mathfrak{su}(6) \times \mathfrak{su}(3)\rangle]$ and the module $R_k[\langle\mathfrak{su}(3)^3\rangle]$ have order 3 twisting; there are fourteen modules with order 2 twist, which we do not list.  The eight centralizers of the affine vertices of the Weyl alcove of $E_8$ are listed in Table~\ref{table-e8}; the corresponding representation modules are all twisted, with order given by the corresponding dual Coxeter label.  Besides these, there are two modules of the form $R_k[\langle\mathfrak{su}(6) \times \mathfrak{su}(3)\rangle]$ and modules $R_k[\langle\mathfrak{su}(3)^3\rangle]$ and $R_k[\langle\mathfrak{su}(3)^3 \times \mathfrak{su}(2)\rangle]$ with order 3 twisting, one module $R_k[\langle\mathfrak{su}(4)^2 \times \mathfrak{su}(2)\rangle]$ with order 4 twisting, and twenty-five modules with order 2 twisting.

\begin{center}
\begin{table}[h]
\caption{Full rank semi-simple centralizers in $E_8$.} \label{table-e8}
\begin{tabular}{|c|l|}
\hline
Vertex & Centralizer \\
\hline
$v_1$  & $Spin(16)/\ZZ/2$ \\
$v_2$  & $SU(9)/\ZZ/3$ \\
$v_3$  & $SU(2) \times_{\ZZ/4} SU(8)$ \\
$v_4$  & $\{SU(2) \times SU(3) \times SU(6)\}/\ZZ/6$ \\
$v_5$  & $SU(5) \times_{\ZZ/5} SU(5)$ \\
$v_6$  & $Spin(10) \times_{\ZZ/4} SU(4)$ \\
$v_7$  & $E_6 \times_{\ZZ/3} SU(3)$ \\
$v_8$  & $E_7 \times_{\ZZ/2} SU(2)$ \\
\hline
\end{tabular}
\end{table}
\end{center}
\nid Here the $\ZZ/2$ action on $Spin(16)$ is the one whose quotient is $sSpin(16)$, and in all cases the finite cyclic group determining the quotient either injects into or surjects onto the center of each factor  of the product of simply connected groups.

\subsection{Abstract presentations of the fusion ring} \label{sec-presentations}

Thanks to Proposition~\ref{prop-free}, the method we used to compute the fusion ring of $G_2$ suffices to determine a presentation of any fusion ring $F_k[G]$ in terms of bases for the various twisted representation modules $R_k[H]$ of centralizers $H$ in $G$.  In particular, there is a level-independent bound on the number of generators of the fusion ideal.

\begin{theorem} \label{thm-finiteness}
For $G$ a compact, simple, simply connected Lie group, there is a positive integer $n_G$ such that for all positive $k$, there is a presentation of the level $k$ fusion ring $F_k[G]$ of the form $F_k[G]=R[G]/(\rho_1(k), \rho_2(k), \ldots, \rho_{n_G}(k))$, where the $\rho_i(k) \in R[G]$ are representations of $G$ depending on the level.  An upper bound on the number of generators $n_G$ of the fusion ideal is $\sum_{i \in \bar{D}} |W_G|/|W_{H_{\widehat{0i}}}|$, where $\bar{D}$ is the nonaffine Dynkin diagram of $G$, and $H_{\widehat{0i}}$ is the centralizer of the edge of the Weyl chamber of $G$ through the fundamental weight $\lambda_i$.
\end{theorem} 

\begin{proof}
By the discussion in section~\ref{sec-resolution}, the fusion ring $F_k[G]$ is the cokernel
$$\coker \left(\bigoplus_{\{i<j\} \in D} R_k[H_{\widehat{\imath\jmath}}] \xra{} \bigoplus_{i \in D} R_k[H_{\hat{\imath}}] \right)$$
This is the homology in degree zero of the complex coming from the twisted Mayer-Vietoris spectral sequence.  The map here is twisted holomorphic induction, described explicitly in terms of Weyl invariants in the weight lattice in Proposition~\ref{prop-resolution}.  The twisted induction map $R_k[H_{\widehat{0ij}}] \xra{} R_k[H_{\widehat{\imath\jmath}}]$ is surjective, where $0$, $i$, and $j$ are distinct nodes of $D$.  Because the complex of twisted representation modules is acyclic in homological degree 1, the cokernel in question is therefore equal to
$$\coker \left(\bigoplus_{0<j \in D} R_k[H_{\widehat{0j}}] \xra{} \bigoplus_{i \in D} R_k[H_{\hat{\imath}}] \right)$$

The map $R_k[H_{\widehat{0j}}] \xra{} R_k[H_{\hat{\jmath}}]$ is surjective.  Let $\{h_1^j, \ldots, h_{b(j)}^j\} \subset R_k[H_{\widehat{0j}}]$ denote a collection of elements mapping to a basis over $R[G]$ for $R_k[H_{\hat{\jmath}}]$, and let $\{g_1^j, \ldots, g_{b(0,j)}^j\} \subset R_k[H_{\widehat{0j}}]$ denote a basis over $R[G]$ of $R_k[H_{\widehat{0j}}]$.  For $h \in R_k[H_{\widehat{0j}}]$, write the differential in the complex as $\partial h = \partial_0 h - \partial_1 h \in R_k[H_{\hat{0}}] \oplus R_k[H_{\hat{\jmath}}]$.  There are representations $c_{rs}^j \in R[G]$ such that $\partial_1 g_r^j = \sum c_{rs}^j \partial_1 h_s^j$.  There results a presentation of the fusion ring:
\begin{align*}
F_k[G] = R[G] / \Big(
& \left(\partial_0 g_1^1 - \sum c_{1s}^1 \partial_0 h_s^1, \partial_0 g_2^1 - \sum c_{2s}^1 \partial_0 h_s^1, \ldots, \partial_0 g_{b(0,1)}^1 - \sum c_{b(0,1)s}^1 \partial_0 h_s^1\right), \\
& \left(\partial_0 g_1^2 - \sum c_{1s}^2 \partial_0 h_s^2, \partial_0 g_2^2 - \sum c_{2s}^2 \partial_0 h_s^2, \ldots, \partial_0 g_{b(0,2)}^2 - \sum c_{b(0,2)s}^2 \partial_0 h_s^2\right), \\ 
& \! \quad \ldots, \\
& \left(\partial_0 g_1^n - \sum c_{1s}^n \partial_0 h_s^n, \partial_0 g_2^n - \sum c_{2s}^n \partial_0 h_s^n, \ldots, \partial_0 g_{b(0,n)}^n - \sum c_{b(0,n)s}^n \partial_0 h_s^n \right) 
\Big)
\end{align*} 
The rank $b(0,j)$ of $R_k[H_{\widehat{0j}}]$ is equal to the ratio $|W_G|/|W_{H_{\widehat{0j}}}|$, where $W_{H_{\widehat{0j}}}$ is the Weyl group of the centralizer of the nonaffine edge $F_{\widehat{0j}}$ of the Weyl alcove.  Note that in this presentation, the differential $\partial_0$ coming from twisted holomorphic induction depends on the level $k$, but the various representations $g_r^j$ and $h_s^j$ only depend on the level modulo the least common multiple of the dual Coxeter labels of the affine Dynkin diagram of $G$.
\end{proof}
\pagebreak

When all of the modules appearing in the complex resolving the fusion ring are untwisted, the number of generators of the fusion ideal can be reduced by making convenient choices of bases.  

\begin{prop}
The ideals defining the fusion rings $F_k[SU(n)]$, $F_k[Sp(n)]$, $F_{2k}[Spin(n)]$, $F_{2k}[G_2]$, $F_{6k}[F_4]$, $F_{6k}[E_6]$, $F_{12k}[E_7]$, and $F_{60k}[E_8]$ are generated by $n_G$ representations, where
$$n_G = \sum_{i \in \bar{D}} \left(\frac{|W_G|}{|W_{H_{\widehat{0i}}}|} - \frac{|W_G|}{|W_{H_{\hat{\imath}}}|}\right)$$
\end{prop}

\begin{proof}
By the discussion in section~\ref{sec-twrepmod}, the indicated levels are those for which all of the modules appearing in the resolution of the fusion ring are untwisted; that is, for those $k$, we have $R_k[H_S] \cong R[H_S]$.  By Steinberg's explicit presentations of representation rings~\cite{steinberg}, there exists a set $\{\mu_1, \ldots, \mu_{b(j)}\} \subset \Lambda_W$ of weights such that the corresponding highest weight irreducible representations $\{\overline{\mu_1}^{H_{\hat{\jmath}}}, \ldots \overline{\mu_{b(j)}}^{H_{\hat{\jmath}}}\} \subset R_k[H_{\hat{\jmath}}]$ form a basis for $R_k[H_{\hat{\jmath}}]$ as an $R[G]$-module.  Under the induction map $\partial_1$ the irreducible representations $\{\overline{\mu_1}^{H_{\widehat{0j}}}, \ldots \overline{\mu_{b(j)}}^{H_{\widehat{0j}}}\} \subset R_k[H_{\widehat{0j}}]$ map to the basis elements $\{\overline{\mu_1}^{H_{\hat{\jmath}}}, \ldots \overline{\mu_{b(j)}}^{H_{\hat{\jmath}}}\} \subset R_k[H_{\hat{\jmath}}]$.  Moreover, Steinberg's construction ensures that the collection $\{\mu_1, \ldots, \mu_{b(j)}\}$ \emph{extends} to a collection $\{\mu_1, \ldots, \mu_{b(j)}, \nu_1, \ldots, \nu_{b(0,j)-b(j)}\}$ of weights whose corresponding highest weight irreducible representations form a basis for $R_k[H_{\widehat{0j}}]$.  In the proof of Theorem~\ref{thm-finiteness}, we may take $h_s^j = \overline{\mu_s}^{H_{\widehat{0j}}}$ and $g_r^j = h_r^j$ for $r \leq b(j)$, so that the submatrix $(c_{rs}^j)_{1 \leq r \leq b(j)}$ of representations is the identity.  As a result, the first $b(j)$ generators of the $b(0,j)$ generators of the fusion ideal induced from $R_k[H_{\widehat{0j}}]$ are zero.
\end{proof}

We expect that the bound from this proposition applies at all levels, not only the levels at which the resolution of the fusion ring is untwisted.  However, even this bound is extremely far from tight---in practice, by exploiting the singular planes for twisted holomorphic induction, one can produce a much smaller set of generators of the fusion ideal.

\bibliography{fusion}

\bibliographystyle{plain}

\end{document}